\documentstyle[art12,leqno]{article}
\title{Non-trivial Linear Systems on Smooth Plane Curves}
\author{Marc Coppens\thanks{This author is related to the University at
Leuven (Celestijnenlaan 200B B-3001 Leuven Belgium) as a Research Fellow.}
 \and Takao Kato\thanks{This author represents his
thanks to Katholieke Universiteit Leuven, Departement Wiskunde for their
kind hospitality during his stay and to the JSPS for financial support.}}
\let\bbb=\bf
\date{}
\def\al{\alpha}
\def\be{\beta}
\def\ga{\gamma}
\def\Ga{{\mit \Gamma}}
\def\ep{\varepsilon}
\def\ze{\zeta}
\def\et{\eta}
\def\ph{\varphi}
\def\Ph{{\mit \Phi}}
\def\ps{\psi}
\def\Ps{{\mit \Psi}}
\def\sus{\subset}
\def\vs#1{\par\vspace{.#1in}}
\def\fracd#1#2{{\displaystyle{#1}\over\displaystyle{#2}}}
\def\Pb{{\bbb P}}
\begin{document}
\newtheorem{lem}{Lemma}[section]
\newtheorem{thm}[lem]{Theorem}
\newtheorem{rem}[lem]{Remark}
\newtheorem{defi}[lem]{Definition}
\newtheorem{prop}[lem]{Proposition}
\newtheorem{cor}[lem]{Corollary}
\newtheorem{const}[lem]{Construction}
\def\theequation{\arabic{section}.\arabic{equation}}
\maketitle
\setcounter{section}{-1}
\section{Introduction}

Let $C$ be a smooth plane curve of degree $d$ defined over an
algebraically closed field $k$.  In \cite{noether}, while studying space
curves, Max Noether considered the following question.  For
$n\in{\bbb Z}_{\ge 1}$ find $\ell (n)\in{\bbb Z}_{\ge 0}$ such that there
exists a linear system $g^{\ell (n)}_n$ on $C$ but no linear system
$g^{\ell (n)+1}_n$ and classify those linear systems $g^{\ell (n)}_n$ on
$C$.  The arguments given by Noether in the answer to this question
contained a gap.  In \cite{cil} C.~Ciliberto gave a complete proof for
Noether's claim using different arguments.  In \cite{harts1}
R.~Hartshorne completed Noether's original arguments by solving the
problem also for integral (not necessarily smooth) plane curves (see
Remark \ref{rem:1}).

The linear systems $g^{\ell (n)}_n$ are either non-special, or special
but not very special, or very special but trivial.  By a very special
(resp. trivial) linear system on a smooth plane curve $C$ we mean:
\vs 2
{\raggedright {\sc Definition}. }
{\it A linear system $g^r_n$ on $C$ is very special if $r\ge 1$ and
$\dim |K_C-g^r_n|\ge 1$. $($Here $K_C$ is a canonical divisor on $C)$.
A base point free complete very special linear system $g^r_n$ on $C$ is
trivial if there exists an integer $m\ge 0$ and an effective divisor $E$
on $C$ of degree $md-n$ such that $g^r_n=|mg^2_d-E|$ and
$r=\fracd{m^2+3m}{2}-(md-n)$.  A complete very special linear system
$g^r_n$ on $C$ is trivial if its associated base point free linear
system is trivial.}
\vs 2
In this paper, we consider the following question.  For
$n\in{\bbb Z}_{\ge 1}$ find $r(n)$ such that there exists a very special
non-trivial complete linear system $g^{r(n)}_n$ on $C$ but no such linear
system $g^{r(n)+1}_n$.  Our main result is the following:
\vs 2
{\raggedright {\sc Theorem}. }
{\it Let $g^r_n$ be a base point free very special non-trivial complete
linear system on $C$. Write $r=\fracd{(x+1)(x+2)}{2}-\be$ with $x, \be$
integers satisfying $x\ge 1, 0\le\be\le x$.  Then
$$
n\ge n(r):=(d-3)(x+3)-\be .
$$}
\vs 2
This theorem only concerns linear systems of dimension $r\ge 2$.  But
 1-dimensional linear systems are studied in \cite{c3}.  From these
results one finds that $C$ has no non-trivial very special linear system
of dimension 1 if $d\le 5$ and for $d\ge 6$, $C$ has non-trivial very
special complete linear systems $g^1_{3d-9}$ but no such linear system
$g^1_n$ with $n<3d-9$.  The proof of this theorem is also effective for
case $r=1$ if one modifies it little bit.

Concerning the original problem one can make the following observation.
For $x\ge d-2$ one has $r>g(C)$ and of course $C$ has no non-trivial
very special linear systems $g^r_n$.  For $x\le d-3$ one has
$(d-3)(x+2)\le (d-3)(x+3)-x$.  So, if the bound $n(r)$ is sharp, then
also the bound $r(n)$ can be found.  Concerning the sharpness of the bound
$n(r)$, we prove it in case ${\rm char}(k)=0$ for $x\le d-6$.  In case
${\rm char}(k)\ne 0$ we prove that there exists smooth plane curves of
degree $d$ with a very special non-trivial $g^r_{n(r)}$ in case
$x\le d-6$.  Finally for the case $x>d-6$ we prove that there exist no
base point free complete very special non-trivial linear systems of
dimension $r$ on $C$.  Hence, at least in case ${\rm char}(k)=0$ the
numbers $r(n)$ are determined.
\vs 2
{\raggedright{\sc Some Notations}}
\vs 1
We write $\Pb_a$ to denote the space of effective divisors of degree $a$
on $\Pb^2$.  If $\Pb$ is a linear subspace of some $\Pb_a$ then we write
$\Pb .C$ for the linear system on $C$ obtained by intersection with
divisors $\Ga\in\Pb$ not containing $C$.  We write $F(\Pb .C)$ for the
fixed point divisor of $\Pb .C$ and $f(\Pb .C)$ for the associated base
point free linear system on $C$, so
$f(\Pb .C)=\{ D-F(\Pb .C):D\in\Pb .C\}$.
If $Z$ is a 0-dimensional subscheme of $\Pb^2$ then $\Pb_a(-Z)$ is the
subspace of divisors $D\in\Pb_a$ with $Z\sus D$.

\setcounter{equation}{0}
\section{Bound on the degree of non-trivial linear systems}

A complete linear system $g^r_n$ on a smooth curve $C$ is called very
special if $r\ge 1$ and $\dim |K_C-g^r_n|\ge 1$.  From now on, $C$ is a
smooth plane curve of degree $d$ and $g^r_n$ is a very special base point
 free linear system on $C$ with $r\ge 2$.
\vs 1
\begin{defi}
$g^r_n$ is called a trivial linear system on $C$ if there exists an
integer $m\ge 0$ and an effective divisor $E$ on $C$ of degree $md-n$
such that $g^r_n=|mg^2_d-E|$ and $r=\fracd{m^2+3m}{2}-(md-n)$.
\label{def:trivial}
\end{defi}
\vs 2
\begin{thm}
Write $r=\fracd{(x+1)(x+2)}{2}-\be$ with $x, \be$ integers satisfying
$x\ge 1, 0\le\be\le x$.  If $g^r_n$ is not trivial, then
$$
n\ge n(r):=(d-3)(x+3)-\be .
$$
\label{thm:1}
\end{thm}
\vs 1
\begin{rem}
{\rm In the proof of this theorem we are going to make use of the main
result of Hartshorne \cite{harts1} which describes the linear systems on
$C$ of maximal dimension with respect to their degrees. The result is as
follows:}

Let $g^r_n$ be a linear system on $C$ $($not necessarily very special$)$.
  Write $g(C)=\fracd{(d-1)(d-2)}{2}$.

{\rm i)} If $n>d(d-3)$ then $r=n-g$ $($the non-special case$)$

{\rm ii)} If $n\le d(d-3)$ then write $n=kd-e$ with
$0\le k\le d-3, 0\le e<d$, one has
\vs 1
$\left\{
\begin{array}{ll}
r\le\fracd{(k-1)(k+2)}{2} & {\rm if\ }e>k+1\\
r\le\fracd{k(k+3)}{2}-e & {\rm if\ }e\le k+1.
\end{array}\right.$
\vs 1
{\rm Hartshorne also gives a description for the case one has equality.
This theorem (a claim originally stated by M. Noether with an incomplete
proof) is also proved in \cite{cil}.  However, Hartshorne also proves the
 theorem for integral plane curves using the concept of generalized
linear systems on Gorenstein curves. We need this more general result in
the proof of Theorem \ref{thm:1}}
\label{rem:1}
\end{rem}
\vs 2
{\it Proof of Theorem\/} \ref{thm:1}.  Assume $g^r_n=rg^1_{n/r}$ and
$n<(x+3)(d-3)-\be$.  Noting $2r=(x+1)(x+2)-2\be\ge x^2+x+2\ge x+3$,
we have $\fracd{(x+3)(d-3)-\be}{r}< 2(d-2)$.  Hence,
$g^1_{n/r}=|g^2_d-P|$ for some $P\in C$.  Since $\dim |rg^1_{n/r}|=r$,
certainly $\dim |2g^1_{n/r}|=2$.  But $\dim|2g^2_d-2P|=3$.
A contradiction.

Since $g^r_n$ is special, there exist an integer $1\le m\le d-3$
and a linear system $\Pb\sus\Pb_m$ such that
$g^r_n=f(\Pb.C)$ and $\Pb$ has no fixed components.  In Lemma
\ref{lem:1} we are going to prove that, because $g^r_n$ is not a multiple
of a pencil, a general element $\Ga$ of $\Pb$ is irreducible.

Now, for each element $\Ga'$ of $\Pb$ we have
$F(\Pb.C)\sus\Ga'$ (inclusion of subschemes of $\Pb^2$).  In
particular $F(\Pb.C)\sus\Ga\cap\Ga'$.  This remark is known
in the literature as Namba's lemma.  As a subscheme of $\Ga$,
$F(\Pb.C)$ is an effective generalized divisor on $\Ga$ (terminology
 of \cite{harts1}).  We find that for each $\Ga'\in\Pb$ with
$\Ga'\ne\Ga$ the residual of $F(\Pb.C)$ in $\Ga\cap\Ga'$ (we denote
it by $\Ga\cap\Ga'-F(\Pb.C)$) is an element of the generalized
complete linear system on $\Ga$ associated to
${\cal O}_{\Ga}(m-F(\Pb.C))$.  Hence, we obtain a generalized linear
 system $g^{r-1}_{m^2-(md-n)}$ on $\Ga$.

Now we are going to apply Hartshorne's theorem (Remark \ref{rem:1}) to
this $g^{r-1}_{m^2-(md-n)}$ on $\Ga$.  Since $g^r_n$ is non-trivial on
$C$, we know that $r>\fracd{m^2+3m}{2}-(md-n)$.  If
$m^2-(md-n)>m(m-3)$, then i) in Remark \ref{rem:1} implies
$r-1\le m^2+n-md-\fracd{(m-1)(m-2)}{2}$ so $r\le\fracd{m^2+3m}{2}-(md-n)$,
 a contradiction.

So $m^2-(md-n)\le m(m-3)$ and we apply ii) in Remark \ref{rem:1}.  We
find $x\le m-3$ and $m^2+n-md\ge mx-\be$, so
$n\ge\ph (m):=-m^2+m(d+x)-\be$.  Since $x+3\le m\le d-3$, we find
$n\ge\ph (x+3)=\ph (d-3)=(d-3)(x+3)-\be=n(r)$.  This completes the proof
of the theorem except for the proof of Lemma \ref{lem:1}.
\vs 2
\begin{lem}
Let $C$ be a smooth plane curve of degree $d$ and let $g^r_n$ be a base
point free complete linear system on $C$ which is not a multiple of a
one-dimensional linear system.  Suppose there exists a linear system
$\Pb\sus\Pb_e$ without fixed component for some $e\le d-1$ such
that $g^r_n=f(\Pb.C)$.  Then the general element of $\Pb$ is irreducible.
\label{lem:1}
\end{lem}
\def\undere{\underline{e}}
\def\underm{\underline{m}}
\vs 2
{\it Proof\/}.  Let $F=F(\Pb .C)=\sum_{j=1}^sn_jP_j$ with $n_j\ge 1$ and
$P_i\ne P\j$ for $i\ne j$.  For $t\in{\bbb Z}_{\ge 1}, \undere =
(e_1,\dots ,e_t)\in ({\bbb Z}_{\ge 1})^t$ with $\sum_{i=1}^te_i=e$ and
$\underm =[m_{ij}]_{1\le i\le t,1\le j\le s}$, let
$$
V(t,\undere ,\underm )=\{ \Ga_1+\cdots +\Ga_t:\Ga_i\in\Pb_{e_i}{\rm\ is\
irreducible\ and\ }i(\Ga_i,C;P_j)=m_{ij}\}.
$$
It is not so difficult to prove that this defines a stratification of
$\Pb_e$ by means of locally closed subspaces.

Since $\Pb$ is irreducible there is a unique triple $(t_0,\undere_0,
\underm_0)$ such that $\Pb\cap V(t_0,\undere_0,\underm_0)$ is an open
non-empty subset of $\Pb$.  In particular, $\Pb\sus
\{ \Ga_1+\cdots +\Ga_{t_0}:\Ga_i\in\Pb_{e_{0i}}{\rm\ and\ }
i(\Ga_i,C;P_j)\ge m_{0ij}\}$.  We need to prove that $t_0=1$,
so assume that $t_0>1$.  Let forget the subscript $0$ from now on.

Let $F_i=\sum_{j=1}^sm_{ij}P_j\sus C$.  For each $D\in g^r_n$ there
exists $\Ga =\Ga_1+\dots +\Ga_t$ with $\Ga_i\in\Pb_i(-F_i)$ and
$D=\Ga .C-(F_1+\dots +F_s)=\sum_{j=1}^t(\Ga_i.C-F_i)$.  Writing
$D_i=\Ga_i.C-F_i\in|e_ig^2_d-F_i|$ we find $D=\sum_{i=1}^tD_i$.
Suppose for some $1\le i\le t$ we have $\dim|e_ig^2_d-F_i|=0$.  If
$\Ga'$ and $\Ga''$ are in $\Pb_i(-F_i)$ then $\Ga'.C=\Ga''.C$, but
$e_i<d$ so $\Ga'=\Ga''$.  This implies that $\Pb_i(-F_i)=\{\Ga_0\}$, but
then $\Ga_0$ is a fixed component of $\Pb$, a contradiction.  Hence, for
$1\le i\le t$, we have $\dim|e_ig^2_d-F_i|\ge 1$.  Now, let $L_i$ be the
irreducible sheaf on $C$ associated to $|e_ig^2_d-F_i|$ and let $L$ be
the irreducible sheaf on $C$ associated to $g^r_n$.  Then
$L=L_1\otimes\cdots\otimes L_t$ and we find that the natural map
$$
H^0(C,L_1)\otimes\cdots\otimes H^0(C,L_t)\to H^0(C,L)
$$
is surjective, while $\dim H^0(C,L_i)\ge 2$ for $1\le i\le t$. From
\cite[Corollary 5.2]{eisen}, it follows that $g^r_n$ is a multiple of a
pencil.  But this is a contradiction.
\vs 2
\begin{rem}
{\rm In \cite{mez} one makes a classification of linear systems on smooth
plane curves for which $r$ is one less than the maximal dimension with
respect to the degree.  In that paper one uses arguments like in
\cite{cil}.  That classification is completely contained in our Theorem
\ref{thm:1}}
\label{rem:2}
\end{rem}
\vs 2
\begin{rem}
{\rm If $r\ge n-\fracd{(d-1)(d-2)}{2}+d-1$ then $\dim |K_C-g^r_n|\ge d-2$.
Hence, $|K_C-g^2_d-g^r_n|=|(d-4)g^2_d-g^r_n|\ne\emptyset$. So in this
case we can assume $m\le d-4$ in the proof of Theorem \ref{thm:1}.
Then, in the proof of Theorem \ref{thm:1}, using Bertini's theorem, we
can prove that, for $D\in g^r_n$ general there exists an irreducible
curve $\Ga$ of degree $d-3$ with $\Ga .C\ge D$ (see argument in
\cite[p.384]{harts1}).  So we don't need the proof of Lemma \ref{lem:1}
under that assumption on $r$.}
\label{rem:25}
\end{rem}
\vs 2
\begin{rem}
{\rm In that case $n=n(r)$, we find $m=x+3$ and $m^2+n-md=mx-\be$. So the
generalized linear system $g^{r-1}_{m^2+n-md}=g^{r-1}_{mx-\be}$ on $\Ga$
is of maximal dimension with respect to its degree.  The description of
those linear systems in \cite{harts1} implies that it is induced by a
family of plane curves of degree $x$ containing some subspace $E\sus\Ga$
of length $\be$.  Writing $Z=F(\Pb .C)\sus\Ga$ we find
$|\Pb_m.\Ga -Z|=|\Pb_x.\Ga -E|$ and so
$Z\in |\Pb_3.\Ga +E|$.  In order to find non-trivial $g^r_{n(r)}$'s
it is interesting to find for a smooth curve $C$ of degree $d$, a curve
$\Ga$ of degree $m$ and a curve $\Ga'$ of degree 3 such that
$\Ga\cap\Ga'\sus C$.  We discuss this in \S 3.  First we solve the
following postulation problem: let $\Ga'$ be the union of 3 distinct
lines $L_1, L_2, L_3$ and let $D_i$ be an effective divisor of degree
$a$ on $L_i$ with $D_i\cap (L_j\cup L_k)=\emptyset$ for
$\{ i,j,k\} =\{ 1,2,3\}$.  Give necessary and sufficient conditions for
the existence of a smooth curve $\Ga$ of degree $a$ such that
$\Ga .L_i=D_i$ for $i=1,2,3$.}
\label{rem:3}
\end{rem}

\setcounter{equation}{0}
\def\binom#1#2{{#1\choose#2}}
\section{Carnot's theorem}

We begin with pointing out the following elementary fact:
\vs 1
\begin{lem}
Let $m(\ge 4)$ and $m_j\ (j=1,\dots ,\ell)$ be positive integers
satisfying $\displaystyle\sum^{\ell}_{j=1}m_j=m$ and let
$\Ph (X)=\displaystyle\sum^m_{j=1}a_jX^j$ be a non-zero polynomial of
degree at most $m$.  If $\Ph (X)$ is divisible by $(X-X_j)^{m_j}$
for $\ell$ distinct values $X_1,\dots ,X_{\ell}$, then the ratio
$a_0:\cdots :a_m$ is uniquely determined.  In particular,
$a_m\ne 0, a_0=(-1)^ma_m\displaystyle\prod^{\ell}_{j=1}X^{m_j}_j$ and
$a_{m-1}=-a_m\displaystyle\sum^{\ell}_{j=1}m_jX_j$.
\label{lem:31}
\end{lem}
\vs 1
Using this, we have the following Carnot's theorem and infinitesimal
Carnot's theorems.  Another generalization of this theorem is given by
Thas et al. \cite{thas} (see also \cite{ver}).  They call this
B.~Segre's generalization of Menelaus' theorem.
\vs 1
\begin{lem}{\rm (Carnot, cf. \cite[Proposition 1.8]{ver},
\cite[p.219]{koe})}
Let $L_1, L_2, L_3$ be lines in $\Pb^2$ so that
$L_1\cap L_2\cap L_3=\emptyset$ and let
$D_i=\displaystyle\sum^{\ell_i}_{j=1}m_{ij}P_{ij},\
(\displaystyle\sum^{\ell_i}_{j=1}m_{ij}=m, i=1,2,3)$ be an effective
divisor on $L_i$ such that $D_{i_1}\cap L_{i_2}=\emptyset$ if
$i_1\ne i_2$.  Assume $(x:y:z)$ is a coordinate system on $\Pb^2$ such
that $L_1, L_2, L_3$ correspond to the coordinate axes $x=0, y=0, z=0$,
respectively.  Let the coordinate of $P_{ij}$ be given by
$(x_{ij}:y_{ij}:z_{ij})$ $($of course $x_{1j}=y_{2j}=z_{3j}=0)$.  Then,
there exists a curve $\Ga$ not containing any one of the lines
$L_1,L_2,L_3$ of degree $m$ satisfying $(\Ga .L_i)=D_i$ $(i=1,2,3)$ if
and only if
\begin{equation}
\prod^{\ell_1}_{j=1}\left(\frac{y_{1j}}{z_{1j}}\right)^{m_{1j}}
\prod^{\ell_2}_{j=1}\left(\frac{z_{2j}}{x_{2j}}\right)^{m_{2j}}
\prod^{\ell_3}_{j=1}\left(\frac{x_{3j}}{y_{3j}}\right)^{m_{3j}}
=(-1)^m.
\label{eq:31}
\end{equation}
\label{lem:car1}
\end{lem}
\vs 1
{\it Proof\/}.  Assume there exists a curve $\Ga$ of degree $m$ not
containing any one of the lines $L_1,L_2,L_3$ satisfying $(\Ga .L_i)=D_i$
$(i=1,2,3)$.  Such a curve is given by
$\Ph (x,y,z)=\displaystyle\sum_{i+j+k=m}a_{ijk}x^iy^jz^k=0$.  In this
description, if $i(\Ga ,L_1;P_{1j})=m_{1j}$, then $\Ph (0,y,z)$ is
divisible by $(y_{1j}z-z_{1j}y)^{m_{1j}}$.
This implies $a_{00m}=(-1)^ma_{0m0}\displaystyle\prod^{\ell_1}_{j=1}
\left(\fracd{y_{1j}}{z_{1j}}\right)^{m_{1j}}$.  Similarly, we have
$a_{m00}=(-1)^ma_{00m}\displaystyle\prod^{\ell_2}_{j=1}
\left(\fracd{z_{2j}}{x_{2j}}\right)^{m_{2j}}$ and
$a_{0m0}=(-1)^ma_{m00}\displaystyle\prod^{\ell_3}_{j=1}
\left(\fracd{x_{3j}}{y_{3j}}\right)^{m_{3j}}$.  Since $\Ga$ does not
contain an intersection point $L_{i_1}\cap L_{i_2}$ for $i_1\ne i_2$, we
have $a_{m00}a_{0m0}a_{00m}\ne 0$.  Hence, we have (\ref{eq:31}).

For the converse, take $a_{m00}=1$.  Then, by (\ref{eq:31}) we can
determine $a_{ijk}$ so that $\Ga$ has the desired property.  This
completes the proof.
\vs 2
Next, we see two infinitesimal cases i.~e. case
$D_{i_1}\cap L_{i_2}\ne\emptyset$ and case
$L_1\cap L_2\cap L_3\ne\emptyset$.
\vs 2
\begin{lem}
Let $L_1, L_2, L_3$ be lines in $\Pb^2$ so that
$L_1\cap L_2\cap L_3=\emptyset$ and let
$D_i=\displaystyle\sum^{\ell_i}_{j=1}m_{ij}P_{ij},
(\displaystyle\sum^{\ell_i}_{j=1}m_{ij}=m)$ be an effective divisor on
$L_i$ such that $m_{11}=m_{21}=1, P_{11}=P_{21}=L_1\cap L_2$ and
$D_i\cap L_3=D_3\cap L_i=\emptyset$ $(i=1,2)$.  Let $(x:y:z)$ and
$(x_{ij}:y_{ij}:z_{ij})$ be as in Lemma {\rm \ref{lem:car1}}.  Then,
there exists a curve $\Ga$, not containing any one of the lines
$L_1, L_2, L_3$, of degree $m$ such that $(\Ga .L_i)=D_i$
$(i=1,2,3)$ and whose tangent line $T$ at $P_{11}=(0:0:1)$ is given by
$\al x+y=0$ $(\al\ne 0)$ if and only if
\begin{equation}
\al
\prod^{\ell_1}_{j=2}\left(\frac{y_{1j}}{z_{1j}}\right)^{m_{1j}}
\prod^{\ell_2}_{j=2}\left(\frac{z_{2j}}{x_{2j}}\right)^{m_{2j}}
\prod^{\ell_3}_{j=1}\left(\frac{x_{3j}}{y_{3j}}\right)^{m_{3j}}
=(-1)^m.
\label{eq:32}
\end{equation}
\label{lem:car2}
\end{lem}
\vs 1
{\it Proof\/}.  We use the same notation as in the proof of Lemma
\ref{lem:car1}.  Assume there exists a curve $\Ga$ having the desired
property. The condition that the tangent line at $P_{11}$ is given by
$\al x+y=0$ implies that $\Ph (0,0,1)=0$ and the linear term of
$\Ph (x,y,1)$ is divisible by $\al x+y$.
Hence, $a_{00m}=0$ and $a_{10,m-1}-\al a_{01,m-1}=0$.  As in the proof
of the previous lemma, we have $a_{m00}a_{0m0}\ne 0, a_{01,m-1}\ne 0$,
$a_{0m0}=(-1)^ma_{m00}\displaystyle\prod^{\ell_3}_{j=1}
\left(\fracd{x_{3j}}{y_{3j}}\right)^{m_{3j}}$,
$a_{01,m-1}=(-1)^{m-1}a_{0m0}\displaystyle\prod^{\ell_1}_{j=2}
\left(\fracd{y_{1j}}{z_{1j}}\right)^{m_{1j}}$ and $\al a_{01,m-1}=
a_{10,m-1}=(-1)^{m-1}a_{m00}\displaystyle\prod^{\ell_2}_{j=2}\left(
\fracd{x_{2j}}{z_{2j}}\right)^{m_{2j}}$. Hence, we have (\ref{eq:32}).

Similar to the proof of the previous lemma, we have the converse.
\vs 2
\begin{lem}
Let $L_1, L_2, L_3$ be lines in $\Pb^2$ so that
$L_1\cap L_2\cap L_3\ne\emptyset$ and let
$D_i=\displaystyle\sum^{\ell_i}_{j=1}m_{ij}P_{ij},
(\displaystyle\sum^{\ell_i}_{j=1}m_{ij}=m)$ be an effective divisor on
$L_i$ such that $D_i\cap L_j=\emptyset$ if $i\ne j$.  Let $(x:y:z)$ be
a coordinate system on $\Pb^2$ such that $L_1, L_2, L_3$ correspond
to the line $y-z=0, y=0, z=0$, respectively.  Let the coordinate
of $P_{ij}$ be given by $(x_{ij}:y_{ij}:z_{ij})$ $($of course
$y_{1j}=z_{1j}, y_{2j}=0, z_{3j}=0)$.  Then,
there exists a curve $\Ga$ of degree $d$ not containing any one of the
lines $L_1,L_2,L_3$ such that $(\Ga .L_i)=D_i$ $(i=1,2,3)$ if and only if
\begin{equation}
\sum^{\ell_1}_{j=1}m_{1j}\frac{x_{1j}}{y_{1j}}-
\sum^{\ell_2}_{j=1}m_{2j}\frac{x_{2j}}{z_{2j}}-
\sum^{\ell_3}_{j=1}m_{3j}\frac{x_{3j}}{y_{3j}}=0.
\label{eq:33}
\end{equation}
\label{lem:car3}
\end{lem}
\vs 1
{\it Proof\/}.  Again, we use the same notation as in the proof of Lemma
\ref{lem:car1}.  Assume there exists a curve $\Ga$ having the desired
property.  Since $\Ga$ does not contain $L_2\cap L_3$, we
have $a_{m00}\ne 0$.  By Lemma \ref{lem:31},
$$
a_{m-1,10}=-\sum^{\ell_3}_{j=1}m_{3j}\frac{x_{3j}}{y_{3j}}a_{m00}\quad
{\rm and}\quad
a_{m-1,01}=-\sum^{\ell_2}_{j=1}m_{2j}\frac{x_{2j}}{z_{2j}}a_{m00}.
$$
To consider the condition on $L_1$, we take the coordinate system
$(\xi :\et :\ze )$ with $\xi =x, \et =y, \ze =z-y$.  Put
$$
\Ps (\xi ,\et ,\ze )=\Ph (\xi ,\et ,\et +\ze )=\sum_{i+j+k=m}a_{ijk}
\xi^i\et^j(\ze +\et )^k=\sum_{i+j+k=m}b_{ijk}\xi^i\et^j\ze^k.
$$
Then, $b_{m00}=a_{m00}$ and $b_{m-1,10}=a_{m-1,10}+a_{m-1,01}$.  In this
description, if $i(\Ga , L_1;P_{1j})=m_{1j}$, then $\Ps (\xi ,\et ,0)$ is
divisible by $(\xi_{1j}\et -\et_{1j}\xi )^{m_{1j}}$.
This implies that $b_{m-1,10}=-\displaystyle\sum^{\ell_1}_{j=1}m_{1j}
\fracd{x_{1j}}{y_{1j}}b_{m00}$.  Then, we have (\ref{eq:33}).

For the converse, noting that if $a_{m00}\ne 0$ then $\Ga$ does not
contains $L_i$ $(i=1,2,3)$ as a component, we can find a $\Ga$
having the desired property.
\vs 2
\begin{rem}
{\rm In each of the lemmas \ref{lem:car1}, \ref{lem:car2} and
\ref{lem:car3}, if (\ref{eq:31}) (resp. (\ref{eq:32}), (\ref{eq:33}))
holds, we can find a smooth curve $\Ga$ of degree $m$ with $\Ga .L_i=D_i$
for $i=1,2,3$.  Indeed, let $\Pb\sus\Pb_m$ be the linear system
of divisors $\Ga$ of degree $m$ on $\Pb^2$ satisfying, as schemes,
$D_i\sus\Ga\cap L_i$.  Clearly $L_1+L_2+L_3+\Pb_{m-3}\sus\Pb$
and we proved that $U=\{\Ga\in\Pb:\Ga{\rm \ does\ not\ contain\ any\
of\ the\ lines\ }L_1, L_2, L_3\}$ is a non-empty open subset of
$\Pb$. Because of Bertini's theorem, for $\Ga\in U$ we have
$L_i\cap\Ga =\{ P_{i1},\dots ,P_{i\ell_i}\}$.
Consider the linear system $\Pb'=\{\Ga\cap\Pb^2\backslash (L_1\cup L_2
\cup L_3):\Ga\in\Pb\}$ on $M=\Pb^2\backslash (L_1\cup L_2\cup L_3)$.
Since $\Ga\cap M\in\Pb'$
for any $\Ga\in\Pb_{m-3}$, $\Pb'$ separates tangent directions and points
on $M$.  Because of Bertini's theorem in arbitrary characteristics (see
\cite{kleiman}), we find that a general element $\Ga$ of $\Pb'$ is
smooth.  So, a general element $\Ga$ of $\Pb$ satisfies ${\rm Sing}(\Ga )
\sus\{ P_{ij}:i=1,2,3{\rm\ and\ }1\le j\le\ell\}$.  But, using
$\Ga'\in\Pb_{m-3}$ suited we find $\Ga =\Ga'+L_1+L_2+L_3$ is smooth at
$P_{ij}$, except for the case $P_{11}=P_{21}$ in Lemma \ref{lem:car2}.
In that case, however, if $\Ga\in U$ then $\Ga$ is smooth at $P_{11}$
because of Bezout's theorem.  This completes the proof of the remark.
(For Bertini's theorem in arbitrary characteristics, one can also use
\cite{greco}.)}
\label{rem:33}
\end{rem}
\setcounter{equation}{0}
\section{Sharpness of the bound}

The next proposition implies that it is enough to solve the postulation
problem mentioned in Remark \ref{rem:3} in order to prove sharpness for
the bound $n(r)$ in Theorem \ref{thm:1}.
\vs 2
\begin{prop}
Let $C$ be a smooth plane curve of degree $d$ and let
$r=\fracd{(x+1)(x+2)}{2}-\be$ with $x, \be\in{\bbb Z}$ satisfying
$4\le x+3\le d-3, 0\le\be\le x$.  Let $n=n(r)=(d-3)(x+3)-\be$.  Suppose
there exists a smooth curve $\Ga$ of degree $m=x+3$, a curve $\Ga'$ of
degree $3$ and an effective divisor $E$ of degree $\be$ on $\Ga$ such
that $Z\sus C$, where the divisor $Z=(\Ga\cap\Ga' )+E$ on $\Ga$ is
considered as a closed subscheme of $\Pb^2$.  Then $|mg^2_d-Z|$ is a
non-trivial $g^r_n$ on $C$.
\label{prop:non-trivial}
\end{prop}
\vs 1
{\it Proof\/}.  We write $E=P_1+\cdots +P_{\be}$. Let $L_1,\dots ,L_{\be}$
be general lines through $P_1,\dots ,P_{\be}$, resp., and let $L_{\be +1}
\dots ,L_x$ be general lines in $\Pb^2$.  If $P\in C$, then we write
$\mu_P(Z)$ for the multiplicity of $Z$ at $P$.
\vs{05}
i) $|mg^2_d-Z|$ is base point free.

Suppose $P$ is a base point for $|mg^2_d-Z|$.  Since
$\Ga .C-Z\in|mg^2_d-Z|$ one finds $P+Z\le\Ga .C$, hence
$i(\Ga ,C;P)>\mu_P(Z)\ge i(\Ga ,\Ga';P)$.  Also
$(\Ga'+\sum_{i=1}^xL_i).C-Z\in|mg^2_d-Z|$, hence
$P\in(\Ga'+\sum_{i=1}^xL_i).C-Z=(\Ga'.C-\Ga'\cap\Ga )+(\sum_{i=1}^xL_i.C
-E)$ (sum of two effective divisors).  Since
$P\not\in\sum_{i=1}^xL_i.C-E$, we find $P\in\Ga'.C-\Ga'\cap\Ga$.  This
implies $i(\Ga',C;P)>i(\Ga ,\Ga';P)$.  But
$i(\Ga ,\Ga';P)\ge\min (i(\Ga ,C;P),i(\Ga',C;P))$ (so called Namba's
lemma), hence we have a contradiction.
\vs{05}
ii) $\dim (|mg^2_d-Z|)\ge r$.

Indeed, $(\Ga'+\Pb_x(-E)).C-Z\sus|mg^2_d-Z|$ and $\dim (\Ga'+\Pb_x(-E))=
\fracd{(x+1)(x+2)}{2}-\be -1$.  But also $\Ga .C-Z\sus|mg^2_d-Z|$ while
$\Ga .C\not\in (\Ga'+\Pb_x(-E)).C$.  This proves the claim.
\vs{05}
iii) $\dim (|mg^2_d-Z|)=r$.

If $\dim (|mg^2_d-Z|)>r$ then on $\Ga$ it induces a linear system
$g^{r'}_{mx-\be}$ with $r'\ge r$.  But Hartshorne's theorem (see 1.3)
implies that this is impossible.

iv) $|mg^2_d-Z|$ is not trivial.

First of all, $|mg^2_d-Z|$ is very special.  Indeed $(d-3-m)g^2_d+Z\sus
|K_C-(mg^2_d-Z)|$.  If $d-3=m$ then from the Riemann-Roch theorem, one
finds $\dim|Z|=1$.

Suppose $|mg^2_d-Z|$ would be trivial, i.~e. $|mg^2_d-Z|=|kg^2_d-F|$ with
$r=\fracd{k^2+3k}{2}-(dk-n)$.  Since $g^r_n$ is very special, one has
$k<d-3$.  Consider $D=(\Ga'.C-\Ga\cap\Ga')+(\sum_{i=1}^xL_i.C-E)$ as in
step i).  There exists $\ga\in\Pb_k(-F)$ with $\ga .C=D+F$.  Because of
Bezout's theorem one has $\ga =\ga'+L_1+\cdots +L_x$.  If $P\in E$ then
$P\not\in\Ga'.C-\Ga'\cap\Ga$, otherwise both $i(\Ga',C;P)>i(\Ga',\Ga ;P)$
and $i(\Ga ,C;P)>i(\Ga',\Ga ;P)$, a contradiction to Namba's lemma.  This
implies $\ga'.C\ge\Ga'.C-\Ga'\cap\Ga$.  Once more from Namba's lemma, we
obtain $\Ga'.\ga'\ge\Ga'.C-\Ga'\cap\Ga$ and so
$$
\deg (\Ga'.\ga')=3(k-x)\ge\deg (\Ga'.C-\Ga'\cap\Ga )=3(d-x-3),
$$
a contradiction to $k<d-3$.
\vs 2
\begin{cor}
Let $C$ be a smooth plane curve of degree $d$.  Assume there exists a
plane curve $\Ga'$ of degree $3$ and a smooth plane curve $\Ga$ of degree
$a\ (4\le a\le d-6)$ such that $\Ga\cap\Ga'\sus C$ $($as schemes$)$.
Then $C$ posseses a non-trivial linear system $g^r_n$ for
$r=\fracd{(a-2)(a-1)}{2}-\be, 0\le\be\le a-3$ and $n=n(r)=a(d-3)-\be$.
\label{cor:1}
\end{cor}
\vs 1
{\it Proof\/}.  $\Ga .C=\Ga\cap\Ga'+D$ for some effective divisor $D$ of
degree $a(d-3)$ on $C$.  But $a(d-3)\ge d-3\ge a+3$, so we can choose an
effective divisor $E\sus D$ of degree $\be$ and then one has to apply
Proposition \ref{prop:non-trivial} to $|ag^2_d-Z|$ for $Z=\Ga\cap\Ga'+E$.
\vs 2
\begin{const}
{\rm Fix $\Ga'\in\Pb_3, \Ga\in\Pb_a\ (a\ge 4)$ general and look at
$\Pb_{a+\ep}(-\Ga\cap\Ga'), (\ep\ge 1)$.  Clearly $\Ga'+\Pb_{a+\ep -3}
\sus\Pb_{a+\ep}(-\Ga\cap\Ga'), \Ga +\Pb_{\ep}\sus\Pb_{a+\ep}
(-\Ga\cap\Ga')$.  Take $P\in\Pb^2\backslash (\Ga\cap\Ga')$.  If
$P\not\in\Ga'$ then using $\Ga'+\Pb_{a+\ep -3}$ one can separate tangent
vectors at $P$.  If $P\not\in\Ga$ then one uses $\Ga +\Pb_{\ep}$.
Kleiman's Bertini theorem \cite{kleiman} in arbitrary characteristics
implies that a general element $C\in\Pb_{a+\ep}(-\Ga\cap\Ga')$ is smooth
outside $\Ga\cap\Ga'$.  But if $\Ga''\in\Pb_{\ep}$ is general then
$\Ga +\Ga''$ is smooth on $\Ga\cap\Ga'$.  This implies that a general
element $C\in\Pb_{a+\ep}(-\Ga\cap\Ga')$ is smooth.  This proves that, for
each $d, 1\le x\le d-6, 0\le\be\le x$, there exists a smooth plane curve
$C$ of degree $d$ possesing a non-trivial $g^r_{n(r)}$ with
$r=\fracd{(x+1)(x+2)}{2}-\be$ and $n(r)=(d-3)(x+3)-\be$.}
\label{const}
\end{const}
\vs 2
In trying to prove this statement for all smooth plane curves of degree
$d$, we only succeeded in case ${\rm char}(k)=0$.  This is the following
theorem.
\vs 2
\begin{thm}
Let $C$ be a smooth plane curve of degree
$d$ over an algebraically closed field of characteristic zero.
Let $d>m\ge 4$.  There exists $\Ga'\in\Pb_3$ and a smooth $\Ga\in\Pb_m$
such that, as schemes,  $\Ga\cap\Ga'\sus C$.
\end{thm}
\vs 1
{\it Proof\/}.  Fix two general lines $L_1, L_2$ in $\Pb^2$, let
$S=L_1\cap L_2$.  We may assume neither $L_1$ nor $L_2$ is a tangent
line of $C$ and $S\not\in C$.  Choose points $P_{11},\dots ,P_{1m}$ on
$C\cap L_1$ and $P_{21},\dots ,P_{2m}$ on $C\cap L_2$.  Choose a general
point $S'$ in $\Pb^2\backslash C\cup L_1\cup L_2$.  The pencil of lines
in $\Pb^2$ through $S'$ induces a base point free $g^1_d$ on $C$.
Because $S'$ is general we have:
\begin{itemize}
\item  If $Q$ is a ramification point of $g^1_d$ then the associated
divisor looks like $2Q+E$ with $Q\not\in E$ and $E$ consists of $d-2$
different points (here we use characteristic zero).

\item  If $Q\in L_i\cap C$ $(i=1,2)$ then $Q$ is not a ramification of
$g^1_d$.  The associated divisor is $Q+E$ with
$E\cap (L_1\cup L_2)=\emptyset$.

\item  The line $SS'$ is not a tangent line of $C$.
\end{itemize}
On the symmetric product $C^{(m)}$ we consider
$V=\{ E\in C^{(m)}:{\rm there\ exists\ }D\in g^1_d\ {\rm with\ }E\le D\}$.
 In terminology of \cite{c2} it is the set $V^1_k(g^1_d)$ and we consider
$V$ with its natural scheme structure.  From Chapter 2 in
{\it loc.~cit.}, it follows that $V$ is a smooth curve.

Let $D_0\in g^1_d$ corresponding to the line $SS'$ and let
$V_0=\{ E\in V:E\le D_0\}$.  We define a map
$\ps :V\backslash V_0\to\Pb^1$ as
follows.  Associated to $E\in V\backslash V_0$ there is a line $L$
through $S'$.  Write $E=P_{31}+\cdots +P_{3m}$.
We distinguish 3 possibilities:
\vs{05}
i)  $E\cap (L_1\cup L_2)=\emptyset$.  Choose coordinates
$x,y,z$ on $\Pb^2$ such that $L_1, L_2, L$ corresponds to the coordinate
axes $x=0, y=0, z=0$, respectively.  Let $(x_{ij}:y_{ij}:z_{ij})$ be the
coordinates of $P_{ij}$ $(i=1,2,3;1\le j\le m)$.  Then
$$
\ps (E)=
\prod^m_{j=1}\left(\frac{y_{1j}}{z_{1j}}\right)
\prod^m_{j=1}\left(\frac{z_{2j}}{x_{2j}}\right)
\prod^m_{j=1}\left(\frac{x_{3j}}{y_{3j}}\right).
$$
As long as we take $L_1, L_2, L$ as axes, this value is independent of
the coordinates.
\vs{05}
ii)  $E\cap (L_1\cup L_2)\ne\emptyset$.  Say
$P_{11}=P_{31}\in E\cap (L_1\cup L_2)$.  Choose coordinates as before
and let $\al x+z=0$ be the equation of the tangent line to $C$ at
$P_{11}$ $(\al\ne 0)$.  Then
$$
\ps (E)=\al
\prod^m_{j=2}\left(\frac{y_{1j}}{z_{1j}}\right)
\prod^m_{j=1}\left(\frac{z_{2j}}{x_{2j}}\right)
\prod^m_{j=2}\left(\frac{x_{3j}}{y_{3j}}\right).
$$
Again taking $L_1, L_2, L$ as axes, this value is independent of
the coordinates.  (Of course this is a function to ${\bbb C}$ and we
consider $\Pb^1={\bbb C}\cup\{\infty\}$.)
\vs{05}
iii)  If $\ps (E)$ is not defined in ${\bbb C}$ then $\ps (E)=\infty$.
\vs{05}
For $E\in V_0$, we define $\ps (E)=(-1)^m$
\vs{05}
This map is a holomorphic map.  Indeed, fixing coordinates $(x:y:z)$ such
that $L_1, L_2$ corresponds to $x=0,y=0$, resp. and $S'=(1:1:0)$, we can
write $z-\ga (x-y)=0$ for the pencil of lines through $S'$ (except for
$SS'$).  If $E\in V\backslash V_0$ and $E$ is a part of a divisor of
$g^1_d$ consisting of $d$ different points, then $\ga$ is a local
coordinate of $V$ at $E$.  In case i) we write down $\ps$ locally as a
holomorphic function in $\ga$.  It is easy to check that $\ps$ is
continuous at $E$ in case ii).

For $E\in V_0$, write $\be z+(x-y)=0$ for the pencil of lines through
$S'$ close to $SS'$.  Let $(x_{3j}:x_{3j}:z_{3j})$ be the coordinates at
the points $P_{3j}$ of $E$.  For $E'\in V$ close to $E$ we have
$E'=\sum_{j=1}^mP'_{3j}$ and coordinates
$(x'_{3j}:\be z'_{3j}+x'_{3j}:z'_{3j})$ at $P'_{3j}$.  Here we can assume
that $x'_{3j}=x'_{3j}(\be ), z'_{3j}=z'_{3j}(\be )$ are holomorphic
functions in $\be$ (local coordinate at $V$ in $E$) and
$x_{3j}=x'_{3j}(0), z_{3j}=z'_{3j}(0)$.

Choose new coordinates $\xi =\be z+x-y, \et =y, \ze =x$.  The coordinates
of $P_{1j}$ are $(0:y_{1j}:\be z_{1j}+y_{1j})$, of $P_{2j}$ are
$(x_{2j}:0:\be z_{2j}+x_{2j})$, of $P'_{3j}$ are
$(x'_{3j}:\be z'_{3j}+x'{3j}:0)$.  We find
\begin{eqnarray}
\ps (E') & = &
\prod^m_{j=1}\frac{y_{1j}}{\be z_{1j}-y_{1j}}
\prod^m_{j=1}\frac{\be z_{2j}+x_{2j}}{x_{2j}}
\prod^m_{j=1}\frac{x'_{3j}}{\be z'_{3j}+x'_{3j}}
\label{eq:v0}
\\
& = & (-1)^m-\left(
\sum^m_{j=1}\frac{z_{3j}}{x_{3j}}-
\sum^m_{j=1}\frac{z_{1j}}{y_{1j}}-
\sum^m_{j=1}\frac{z_{2j}}{x_{2j}}\right)\be
+o(\be ).\nonumber
\end{eqnarray}
Hence, $\ps$ is continuous at $E$.

Since $V$ is smooth and $\ps$ is continuous on $V$ and holomorphic except
for a finite number of points, $\ps$ is a holomorphic map $V\to\Pb^1$.

At some component of $V$, $\ps$ is not constant.
Indeed, look at a fibre $2Q+E\in g^1_d$ with $E\in C^{(d-2)}$.  Take a
close fibre $P_1+P_2+E'$ with $P_1, P_2$ close to $Q$.  Choose $F\le E'$
with $\deg(F)=m-1$ and consider $P_1+F\in V$.  Let $W$ be the irreducible
component of $V$ containing $P_1+F$.  Using monodromy one finds
$P_2+F\in W$.  But clearly $\ps (P_1+F)\ne\ps (P_2+F)$, hence
$\ps :W\to\Pb^1$ is a covering.
In particular $\ps^{-1}((-1)^m)\ne\emptyset$.

If for some $E\in W\backslash V_0$ we have $\ps (E)=(-1)^m$ then the
theorem follows from Lemmas \ref{lem:car1}, \ref{lem:car2} and Remark
\ref{rem:33}.  So, we have to take a closer look to $\ps$ at $V_0$.
By the equation (\ref{eq:v0}), if $E\in V_0$ is not a simple zero of
$\ps -(-1)^m$ then the theorem follows from Lemma \ref{lem:car3} and
Remark \ref{rem:33}.

Suppose that each zero of $\ps -(-1)^m$ belonging to $V_0$ is simple.
Then $\ps -(-1)^m$ has exactly $\binom{d}{m}$ zeros at those points.
Now we look at zeros of $\ps$ on $V\backslash V_0$.  The number of zeros
is finite.  For case i) there is none.  For case ii) we have two
possiblities.  If $E\in V$ corresponds to a line $L$ through $S'$
containing one of the points $P_{21},\dots ,P_{2m}$ but
$E\cap (L_1\cup L_2)=\emptyset$.  There are $m\binom{d-1}{m}$ such
possibilities.  If $E\in V$ corresponds to a line $L$ through $S'$ not
containing any of the points $P_{11},\dots ,P_{1m}$ but
$E\cap L_1\ne\emptyset$.  There are $(d-m)\binom{d-1}{m-1}$ such
possibilities.  So, on the components of $V$ where $\ps$ is not constant,
$\ps$ has at least $m\binom{d-1}{m}+(d-m)\binom{d-1}{m-1}$ zeros.
But this number is greater than $\binom{d}{m}$, so $\ps -(-1)^m$ has a
zero on $V\backslash V_0$. This completes the proof of the theorem.
\vs 2
\begin{rem}
{\rm In order to obtain the bound $r(n)$ mentioned in the introduction,
we have to prove that $C$ possesses no base point free very special
non-trivial linear systems $g^r_n$ with $r\ge\fracd{(d-4)(d-3)}{2}-(d-5)$
(i.~e. $x\ge d-5$).  In the introduction we already noticed that
$x\le d-3$.  Assume $g^r_n$ is a very special non-trivial linear system.
{}From Theorem \ref{thm:1} we find $n\ge n((d-4)(d-5)/2-(d-5))=d^2-6d+11$.
But then $\deg (K_C-g^r_n)\le (d-1)(d-2)-2-(d^2-6d+11)=3d-11$.  However,
very special linear systems $g^s_m$ of degree $m\le 3d-11$ are trivial.
So, the associated base point free linear system $g^s_m$ of $|K_C-g^r_n|$
is of type $|ag^2_d-E|$ with $a\le d-4$, $E$ effective and
$\dim|K_C-g^r_n|=\fracd{a^2+3a}{2}-\deg E$.  If $E\ne\emptyset$, then for
$P\in E$ one has $\dim|ag^2_d-E+P|>\dim|ag^2_d-E|$, so $\dim|g^r_n-P|=r$.
This implies that $g^r_n$ has a base point, hence $E=\emptyset$.  But
then, $g^r_n=|(d-3-a)g^2_d-F|$ and since $\dim|ag^2_d+F|=\dim|ag^2_d|$,
we have $r=\fracd{(d-3-a)^2+3(d-3-a)}{2}-\deg F$.  This is a
contradiction to the fact that $g^r_n$ would be non-trivial.}
\end{rem}
\vs 1
\begin{rem}
{\rm It would be interesting to find an answer to the following questions:
For which values of $n$ do there exist non-trivial base point free very
special linear systems $g^r_n$ on a (general) smooth plane curve.
Classify those linear systems and study $W^r_n$ on $J(C)$.  More
concretely, is the subscheme of $W^r_{n(r)}$ corresponding to non-trivial
linear systems irreducible ?  What are the dimension of those irreducible
components ?  And so on.
}
\end{rem}

\begin{flushleft}
Katholieke Industri\"{e}le Hogeschool der Kempen\\
Campus H.~I.~Kempen Kleinhoefstraat 4\\
B 2440 Geel, Belgium\\
\vspace{.1in}
Department of Mathematics, Faculty of Science,\\
Yamaguchi University, Yamaguchi, 753 Japan
\end{flushleft}

\begin{thebibliography}{99}
\hyphenation{auto-mor-ph-ism}

\bibitem{cil} Ciliberto,~C., Alcune applicazioni di un classico
procedimento di Castelnuovo. Seminari di Geometria 1982--83,
Universit\`{a} di Bologna, 17--43.

\bibitem{c2} Coppens, M., An infinitesimal study of secant space divisors.
 preprint

\bibitem{c3} Coppens, M., The existence of base point free linear systems
 on smooth plane curves. preprint

\bibitem{eisen} Eisenbud, D., Linear sections of determinantal varieties.
 Amer. J. Math. {\bf 110}, 541--575 (1988)

\bibitem{greco}  Greco, S. and Valabrega, P., On the singular locus of a
general complete intersection through a variety in projective space.
 Boll. U. M. I. Algebra e Geometria, {\bf 2}, 113--145 (1983)

\bibitem{harts1} Hartshorne, R., Generalized divisors on Gorenstein curves
and a theorem of Noether. J. Math. Kyoto Univ. {\bf 26},
 375--386 (1986)

\bibitem{kleiman}  Kleiman, S., The transversality of a general translate.
 Compos. Math. {\bf 28}, 287--297 (1974)

\bibitem{koe} K\"{o}tter, E., Die Entwickelung der synthetischen
Geometrie. Jahresbericht der D.M.V. {\bf 5}, 1--481 (1896).

\bibitem{mez} Mezzetti, E., Differential-geometric methods for the lifting
problem and linear systems on plane curves. preprint

\bibitem{noether} Noether, M., Zur Grundlegung der Theorie der
algebraischen Raumcurven.  Abh. K\"{o}niglichen Preus. Akad. der
Wissenschaften, 271--318 (1883)

\bibitem{thas} Thas, J.~A., Cameron, P.~J. and Blokhuis, A., On a
generalization of a theorem of B.~Segre. Geometriae Dedicata
{\bf 43}, 299--305 (1992)

\bibitem{ver} Vermeulen, A.~M., Weierstrass points of weight two on curves
of genus three. Ph-D thesis Amsterdam (1983).

\end{thebibliography}
\end{document}